\documentclass[11pt]{article}

\usepackage{amscd,amstex,amssymb}

\def\eqref#1{(\ref{#1})}

\newcommand{\arrow}{{\:\longrightarrow\:}}
\newcommand{\Z}{{\Bbb Z}}
\newcommand{\C}{{\Bbb C}}
\newcommand{\R}{{\Bbb R}}

\newcommand{\Q}{{\Bbb Q}}

\newcommand{\1}{\sqrt{-1}\:}
\newcommand{\inangles}[1]{{\langle #1\rangle}}

\renewcommand{\c}[1]{{\cal #1}}

\renewcommand{\phi}{\varphi}
\renewcommand{\epsilon}{\varepsilon}

\newcommand{\comment}[1]{{}}

\def\blacksquare{\hbox{\vrule width 4pt height 4pt depth 0pt}}
\def\endproof{\blacksquare}

\makeatletter

\@ifundefined{Bbb}
     {\newcommand{\Bbb}[1]{{\mathbb #1}}}%
{}%     {\edef\Bbb#1{{\Bbb #1}}}

\newcommand{\ps@verbit}{%
  \renewcommand{\@oddhead}{%
          \scriptsize
          {Algebraic structures on hyperk\"ahler manifolds}
          \hfil\tiny {September 1996}}
  \renewcommand{\@evenhead}{\@oddhead}
  \renewcommand{\@oddfoot}{\hfil\thepage\hfil}
  \renewcommand{\@evenfoot}{\@oddfoot}}
 
\newcounter{Mycounter}[section]
\newcounter{lemma}[section]
\renewcommand{\thelemma}{{Lemma \thesection.\arabic{lemma}}}
\newcommand{\lemma}{%
     \setcounter{lemma}{\value{Mycounter}}
     \refstepcounter{lemma}
     \stepcounter{Mycounter}
     {\bf \thelemma:\ }}

\newcounter{claim}[section]
\newcounter{sublemma}[section]
\newcounter{corollary}[section]
\renewcommand{\thecorollary}{{Corollary \thesection.\arabic{corollary}}}
\newcommand{\corollary}{%
     \setcounter{corollary}{\value{Mycounter}}
     \refstepcounter{corollary}
     \stepcounter{Mycounter}
     {\bf \thecorollary:\ }}

\newcounter{theorem}[section]
\renewcommand{\thetheorem}{{Theorem \thesection.\arabic{theorem}}}
\newcommand{\theorem}{%
     \setcounter{theorem}{\value{Mycounter}}
     \refstepcounter{theorem}
     \stepcounter{Mycounter}
     {\bf \thetheorem:\ }}

\newcounter{proposition}[section]
\renewcommand{\theproposition}
       {{Proposition \thesection.\arabic{proposition}}}
\newcommand{\proposition}{%
     \setcounter{proposition}{\value{Mycounter}}
     \refstepcounter{proposition}
     \stepcounter{Mycounter}
     {\bf \theproposition:\ }}

\newcounter{definition}[section]
\renewcommand{\thedefinition}
       {{Definition \thesection.\arabic{definition}}}
\newcommand{\definition}{%
     \setcounter{definition}{\value{Mycounter}}
     \refstepcounter{definition}
     \stepcounter{Mycounter}
     {\bf \thedefinition:\ }}

\newcounter{example}[section]
\newcounter{remark}[section]
\renewcommand{\theremark}{{Remark \thesection.\arabic{remark}}}
\newcommand{\remark}{%
     \setcounter{remark}{\value{Mycounter}}
     \refstepcounter{remark}
     \stepcounter{Mycounter}
     {\bf \theremark:\ }}

\newcounter{problem}[section]
\newcounter{question}[section]
\@addtoreset{equation}{section}
\@addtoreset{footnote}{section}
\makeatother

\begin{document}

\begin{center}
{\Large\bf
	Algebraic structures on hyperk\"ahler 
manifolds.}\\[4mm]
Misha Verbitsky,\\[4mm]
{\tt verbit@@thelema.dnttm.rssi.ru}
\end{center}

\hfill

{\small 
\hspace{0.2\linewidth}
\begin{minipage}[t]{0.7\linewidth}
Let $M$ be a compact hyperk\"ahler manifold. The hyperk\"ahler
structure equips $M$ with a set $\c R$ of complex structures 
parametrized by $\C P^1$, called 
{\bf the set of induced complex structures.} 
It was known previously that induced complex structures are
{\it non-algebraic}, except may be a countable set. We prove that
a {\it countable} set of induced complex structures is
{\it algebraic}, and this set is dense in $\c R$. 
A more general version of this theorem was proven by A. Fujiki.
\end{minipage}
}

\hfill

\hfill

The structure of this paper is following. In Subsection
\ref{_hype_defi_Subsection_}
we define hyperk\"ahler manifolds and induced complex structures.
In Subsection \ref{_gene_appe_Subsection_}, we define
induced complex structures of general type.
The main result of this paper and its proof are
given in Section \ref{_algebra_Section_}.

%%%%%%%%%%%%%%%%%%%%%%%%%%%%%%%%%%%%%%%%

\section{Introduction.}
\label{_Intro_Section_}

%%%%%%%%%%%%%%%%%%%%%%%%%%%%%%%%%%%%%%%%

We give the basic definitions and cite the results relevant to 
this paper.

%%%%%%%%%%%%%%%%%%%%%%%%%%%%%%%%%%%%%%%%%%%%%%%%

\subsection{Hyperk\"ahler manifolds}
\label{_hype_defi_Subsection_}

%%%%%%%%%%%%%%%%%%%%%%%%%%%%%%%%%%%%%%%%%%%%%%%%

\hfill

%%%%%%%%%%%%%%%%%%%%%%%%%%%%%%%%%%%%%%%%%%%%%%%%%%%%%%%%%%%%%%%%%
\definition \label{_hyperkahler_manifold_Definition_} %%%%%%%%%%
(\cite{_Besse:Einst_Manifo_}) A {\bf hyperk\"ahler manifold} is a
Riemannian manifold $M$ endowed with three complex structures $I$, $J$
and $K$, such that the following holds.
 
\begin{description}
\item[(i)]  the metric on $M$ is K\"ahler with respect to these complex 
structures and
 
\item[(ii)] $I$, $J$ and $K$, considered as  endomorphisms
of a real tangent bundle, satisfy the relation 
$I\circ J=-J\circ I = K$.
\end{description}

\hfill 

The notion of a hyperk\"ahler manifold was 
introduced by E. Calabi (\cite{_Calabi_}).

\hfill
 
Clearly, a hyperk\"ahler manifold has the natural action of
quaternion algebra ${\Bbb H}$ in its real tangent bundle $TM$. 
Therefore its complex dimension is even.
For each quaternion $L\in \Bbb H$, $L^2=-1$,
the corresponding automorphism of $TM$ is an almost complex
structure. It is easy to check that this almost 
complex structure is integrable (\cite{_Besse:Einst_Manifo_}).

\hfill

%%%%%%%%%%%%%%%%%%%%%%%%%%%%%%%%%%%%%%%%%%%%%%%%%%%%%%%%%%%%
\definition \label{_indu_comple_str_Definition_} 
Let $M$ be a hyperk\"ahler manifold, $L$ a quaternion satisfying
$L^2=-1$. The corresponding complex structure on $M$ is called
{\bf an induced complex structure}. The $M$ considered as a complex
manifold is denoted by $(M, L)$.

\hfill

Let $M$ be a hyperk\"ahler manifold. We identify the group $SU(2)$
with the group of unitary quaternions. This gives a canonical 
action of $SU(2)$ on the tangent bundle, and all its tensor
powers. In particular, we obtain a natural action of $SU(2)$
on the bundle of differential forms. 

\hfill

%%%%%%%%%%%%%%%%%%%%%%%%%%%%%%%%%%%%%%%%%%%%%%%%%%%%%%%%%%%%
\lemma \label{_SU(2)_commu_Laplace_Lemma_}
The action of $SU(2)$ on differential forms commutes
with the Laplacian.

{\bf Proof:} This is Proposition 1.1
of \cite{Verbitsky:Symplectic_II_}. \endproof

Thus, for compact $M$, we may speak of the natural action of
$SU(2)$ in cohomology.

%%%%%%%%%%%%%%%%%%%%%%%%%%%%%%%%%%%%%%%%%%%%%%%%%%%%%%%%%%%%

\subsection{Appendix: Induced complex structures of general type.}
\label{_gene_appe_Subsection_}

%%%%%%%%%%%%%%%%%%%%%%%%%%%%%%%%%%%%%%%%%%%%%%%%%%%%%%%%%%%%

We cite results useful for understanding
of behaviour of induced complex structures. The function
of this appendix
is illustrative. We do not refer to 
this subsection in the body of the article.
This appendix is perfectly safe to skip.

\hfill

%%%%%%%%%%%%%%%%%%%%%%%%%%%%%%%%%%%%%%%%%%%%%%%%%%%%%%%%%%%%
\definition \label{_generic_manifolds_Definition_} %%%%%%%%%
Let $M$ be a hyperk\"ahler manifold, $\c R$ the set of
all induced complex structures. With each $I\in \c R$, we
associate the Hodge decomposition $H^*(M) = \oplus H^{p,q}_I(M)$
on the cohomology of $M$.
We say that $I$ is {\bf 
of general type} when all elements of the group
\[ H^{p,p}_I(M)\cap H^{2p}(M,\Z)\] 
are $SU(2)$-invariant.

\hfill

As \ref{_gene_type_co_div_by2_Remark_} below
implies, the manifolds $(M, I)$ have no Weil divisors
when $I$ is of general type.
In particular, induced complex structures of 
general type are never algebraic.

\hfill

%%%%%%%%%%%%%%%%%%%%%%%%%%%%%%%%%%%%%%%%%%%%%%%%%%%%%%%%%%%%
\proposition \label{_generic_are_dense_Proposition_} %%%%%%%
Let $M$ be a hyperk\"ahler manifold and $\c R$
be the set of induced complex structures over $M$. 
Let $\c R_{ng}\subset \c R$ be the set of all 
induced complex structures {\bf not} of general type.
Then $\c R_{ng}$ is no more than countable.

{\bf Proof:} This is Proposition 2.2 from
\cite{Verbitsky:Symplectic_II_}
\endproof

\hfill

Let $M$ be a compact hyperk\"ahler manifold, $dim_\R M =2m$.

\hfill

%%%%%%%%%%%%%%%%%%%%%%%%%%%%%%%%%%%%%%%%%%%%%%%%%
\definition\label{_trianalytic_Definition_} %%%%%%%%%%
Let $N\subset M$ be a closed subset of $M$. Then $N$ is
called {\bf trianalytic} if $N$ is an analytic subset 
of $(M,L)$ for any induced complex structure $L$.

\hfill

Let $I$ be an induced complex structure on $M$,
and $N\subset(M,I)$ be
a closed analytic subvariety of $(M,I)$, $dim_\C N= n$.
Denote by $[N]\in H_{2n}(M)$ the homology class 
represented by $N$. Let $\inangles N\in H^{2m-2n}(M)$ denote 
the Poincare dual cohomology class. Recall that
the hyperk\"ahler structure induces the action of 
the group $SU(2)$ on the space $H^{2m-2n}(M)$.

\hfill

%%%%%%%%%%%%%%%%%%%%%%%%%%%%%%%%%%%%%%%%%%%%%%%%%%%%%%%%%%%%%%%%%
\theorem\label{_G_M_invariant_implies_trianalytic_Theorem_} %%%%%
Assume that $\inangles N\in  H^{2m-2n}(M)$ is invariant with respect
to the action of $SU(2)$ on $H^{2m-2n}(M)$. Then $N$ is trianalytic.

{\bf Proof:} This is Theorem 4.1 of 
\cite{Verbitsky:Symplectic_II_}.
\endproof

%%%%%%%%%%%%%%%%%%%%%%%%%%%%%%%%%%%%%%%%%%%%%%%%
\remark \label{_triana_dim_div_4_Remark_}
Trianalytic subvarieties have an action of quaternion algebra in
the tangent bundle. In particular,
the real dimension of such subvarieties is divisible by 4.

\hfill

\ref{_G_M_invariant_implies_trianalytic_Theorem_} has the following
immediate corollary, also proven by a Fujiki
(\cite{_Fujiki_}, Theorem 4.8 (1)):

%%%%%%%%%%%%%%%%%%%%%%%%%%%%%%%%%%%%%%%%%%%%%%%%%%%%%%%%%%%%%
\corollary \label{_hyperkae_embeddings_Corollary_} %%%%%%%%%%
Let $M$ be a compact hyperk\"ahler manifold,
$I$ induced complex structure 
of general type, and $S\subset (M,I)$ its complex analytic
subvariety. Then $S$ is trianalytic.

\endproof

%%%%%%%%%%%%%%%%%%%%%%%%%%%%%%%%%%%%%%%%%%%%%%%%%%%%%%%%%%%%
\remark \label{_gene_type_co_div_by2_Remark_}
{}From \ref{_hyperkae_embeddings_Corollary_} and
\ref{_triana_dim_div_4_Remark_}, it follows that
a holomorphically symplectic manifold of general type
has no closed complex analytic subvarieties of odd dimension;
in particular, such a manifold has no divisors.

%%%%%%%%%%%%%%%%%%%%%%%%%%%%%%%%%%%%%%%%%%%%%%%%%%%%%%%%%%%%

\section{Induced complex structures which are algebraic}
\label{_algebra_Section_}

%%%%%%%%%%%%%%%%%%%%%%%%%%%%%%%%%%%%%%%%%%%%%%%%%%%%%%%%%%%%

\hfill

%%%%%%%%%%%%%%%%%%%%%%%%%%%%%%%%%%%%
\definition
Let $M$ be a compact hyperk\"ahler manifold. Then $M$ is called
{\bf simple} if $M$ is simply connected and cannot be 
non-trivially represented as
a direct product of hyperk\"ahler manifolds.

\hfill

The general version of the 
following theorem was proven by
A. Fujiki (\cite{_Fujiki_}, Theorem 4.8 (2)).
Let $\pi:\; \c M \arrow S$ be a deformation of a simple holomorphically
simplectic manifold, with arbitrary base of positive dimension.
Assume that $\c M\arrow S$ is not isotrivial (not trivial on periods). 
Fujiki proves that for a dense subset $S_a\subset S$, the fibers 
$\pi^{-1}(s_a)$ are algebraic for all $s_a\in S_a$. 
I am grateful to Daniel Huybrechts, who provided me with this
reference. Also, a similar (but weaker) result was 
proven by F. Campana (\cite{_Campana_}).

\hfill

%%%%%%%%%%%%%%%%%%%%%%%%%%%%%%%%%%%%%%%%%%%%%%%%
\theorem \label{_alge_dense_Theorem_}
Let $M$ be a compact simple hyperk\"ahler manifold and
$\c R$ be the set of induced complex structures 
$\c R \cong \C P^1$. Let $\c R_{alg}\subset \c R$
be the set of all $L\in \c R$ such that the complex manifold
$(M, L)$ is algebraic. Then $\c R_{alg}$ is countable and dense in
$\c R$. 

\hfill

{\bf Proof:} For each $L\in R$, consider the K\"ahler cone of
$(M,L)$, denoted by $K(L)\subset H^2(M, \R)$. By definition,
$K(L)$ is the set of all cohomology 
classes $\omega\subset H^2(M,\R)$ which are K\"ahler classes
with respect to some metric on $(M,L)$. Let
\[
   K:= \bigcup\limits_{L\in \c R} K(L).
\]
By \cite{_main_}, Lemma 5.6, $K$ is open in $H^2(M, \R)$.
Therefore, the intersection $H^2(M, \Q) \cap K$ is dense in 
$K$. By Kodaira, a compact K\"ahler manifold is algebraic
if and only if there exist a rational K\"ahler class on $M$
(\cite{_Griffi_Harri_}).
By \ref{_K_proje_on_R_Lemma_} below, every cohomology class
$\omega \in K$ corresponds to a unique 
induced complex structure $I(\omega)\in \c R$ such that
$\omega$ is K\"ahler with respect to $I(\omega)$. We also prove that
thus obtained map $\pi:\; K \arrow \c R$ is continuous. 
Thus, $\pi\left( H^2(M, \Q \cap K)\right)$ is dense in $\c R$. On
the other hand, by Kodaira, $\pi\left( H^2(M, \Q \cap K)\right)$ 
coinsides with $\c R_{alg}$. This proves 
\ref{_alge_dense_Theorem_}.
\endproof

\hfill

The following lemma is implicit from \cite{_main_}. We decided
to spell out its proof, for clarity; for missing details the
reader is referred to \cite{_main_}.

\hfill

%%%%%%%%%%%%%%%%%%%%%%%%%%%%%%%%%%%%%%%%%%%%%%%%
\lemma \label{_K_proje_on_R_Lemma_}
Let $\omega\in K$ be a cohomology class which is K\"ahler
with respect to $I\in \c R$. Then 
\begin{description} 
\item [(i)]
   such $I$ is unique,
\item[(ii)]
and the obtained map $\pi:\; K \arrow \c R$ is continuous.
\end{description}

\hfill

{\bf Proof:} Consider the positively defined scalar 
product on the cohomology space $H^2(M)$ induced by the 
standard scalar product on harmonic forms. This scalar
product is clearly $SU(2)$-invariant. 
For an induced complex structure $L$, denote by 
$\omega_L\in H^2(M, \R)$ the
K\"ahler class of the K\"ahler structure associated with $L$
and the hyperk\"ahler structure. 
The corresponding harmonic form can be expressed as 
$\omega_L(x, y) = (x, L(y))$, where $(\cdot,\cdot)$ is the
Riemannian form on $M$.
Let $V\subset H^2(M, \R)$ be the 
3-dimensional subspace generated by $\omega_L$, for all $L\in \c R$
(see \cite{_main_}, Section 4), and $p:\; H^2(M) \arrow V$ be the
orthogonal projection to $V$. For a K\"ahler class
$\omega$ on $(M, I)$, the product $(\omega, \omega_I)$ 
is positive, by Hodge--Riemann relations 
(\cite{_Griffi_Harri_}). Thus, for all $\omega\in K$,
the vector $p(\omega)\in V$ is non-zero. Now,
for each $I\in \c R$, the intersection $H^{1,1}_I(M) \cap V$
is one-dimensional, because in notation of 
\ref{_hyperkahler_manifold_Definition_}, the cohomology class 
$\omega_J+\1 \omega_K\in V$ belongs to $H^{2,0}_I(M)$
and $\omega_J-\1 \omega_K\in V$ belongs to $H^{0,2}_I(M)$.
Therefore, $\omega$ and $\omega_I$ are collinear. 
Clearly, a complex structure $L$ is uniquely defined 
by the form $\omega_L$, and
for all $L\in \c R$, the vectors $\omega_L\in V$ all
have the same length. Thus, for each line $l\in V$,
there exist no more than two induced complex structures
$L\in \c R$ satisfying $l\in H^{1,1}_L(M)$. It is easy
to check that these induced complex structures are opposite:
we have $H^{1,1}_{L}(M) = H^{1,1}_{-L}(M)$.
This implies
that in $\c R$, only $L=I$ and $L=-I$ satisfy 
$\omega\in H^{1,1}_L(M)$. On the other hand,
$\omega_{L} = - \omega_{-L}$. Thus, 
the numbers $(\omega, \omega_I)$ and $(\omega, \omega_{-I})$
cannot be positive simultaneously. This implies that
$\omega$ cannot be a K\"ahler class for $I$ and $-I$ 
at the same time. We proved \ref{_K_proje_on_R_Lemma_} (i).

\hfill

To prove \ref{_K_proje_on_R_Lemma_} (ii), consider the
composition $s$ of $p:\; K \arrow V\backslash 0$ and
the natural projection map from $V\backslash 0 $ to the sphere 
$S^2\subset V \cong \R^3$. If we identify $S^2$ with 
$\c R\cong\C P^1$, we find that $\pi$ is equal
to $s$ (see, for instance, \cite{_main_}, the proof
of Sublemma 5.6). On the other hand, $s$ is 
continous by construction. This proves
\ref{_K_proje_on_R_Lemma_} (ii). \endproof

\hfill

I am grateful to Daniel Huybrechts, who provided me with reference
to \cite{_Campana_} and \cite{_Fujiki_} after the preliminary
version of this paper appeared in alg-geom preprint server.

\hfill

\end{document}